\documentclass[10pt,oneside,leqno]{amsart}
\def\C{{\mathbb C}}
\def\N{{\mathbb N}}
\def\Q{{\mathbb Q}}
\def\Z{{\mathbb Z}}
\def\R{{\mathbb R}}
\def\D{{\mathbb D}}
\def\A{{\mathbb A}}
\def\NQ{{\mathbb{R}\setminus\mathbb{Q}}}
\usepackage{geometry}
\usepackage{enumerate}
\theoremstyle{theorem}\newtheorem{defin}{Definition} \theoremstyle{theorem}\newtheorem{lemma}[defin]{Lemma}
\theoremstyle{theorem}\newtheorem{cor}[defin]{Corollary} \theoremstyle{theorem}\newtheorem{theo}[defin]{Theorem}
\theoremstyle{theorem} \theoremstyle{remark}\newtheorem{remark}[defin]{Remark}
\def\prop{\operatorname{Prop}}
\def\aut{\operatorname{Aut}}
\def\ind{\operatorname{Ind}}
\title[Proper holomorphic mappings]
{Proper holomorphic mappings in the special class of Reinhardt domains}
\author[\L.~Kosi\'nski]
{\L ukasz Kosi\'nski}
\address{Instytut Matematyki, Uniwersytet Jagiello\'nski, Reymonta 4,
30-059\\ Krak\'ow, Poland}\email{lukasz.kosinski@im.uj.edu.pl} 

\begin{document}
\baselineskip=16pt
\begin{abstract}
\baselineskip=16pt A complete characterization of proper holomorphic mappings between domains from the class of all pseudoconvex Reinhardt
domains in $\C^2$ with the logarithmic image equal to a strip or a half-plane is given.
\end{abstract}
\maketitle
\section{Statement of results}
We adopt here the standard notations from complex analysis. Given $\gamma=(\gamma_1,\gamma_2)\in\R^2$ and $z=(z_1,z_2)\in\C^2$ for which it
makes sense we put $|z^{\gamma}|=|z_1|^{\gamma_1}|z_2|^{\gamma_2}.$ The unit disc in $\C$ is denoted by $\D$ and the set of proper holomorphic
mappings between domains $D,G\subset\C^n$ is denoted by $\prop(D,G).$

In this paper we deal with the pseudoconvex Reinhardt domains in $\C^2$ whose logarithmic image is equal to a strip or a half-plane. Observe
that such domains are always algebraically equivalent to domains of the form
\begin{align*}D_{\alpha,r^-,r^+}:=\{z\in\C^2:\
r^-<|z^{\alpha}|<r^+\},\quad \end{align*} where $\alpha=(\alpha_1,\alpha_2)\in(\R^2)_*,\ 0<r^+<\infty,\ -\infty<r^-<r^+.$

We say that $D_{\alpha,r^-,r^+}$ is of the \textit{irrational type} if $\alpha_1/\alpha_2\in\NQ.$ In the other case $D_{\alpha,r^-,r^+}$ is said
to be of the \textit{rational type}.

Recall that if $r^-<0<r^+,\ \alpha\in(\R^2)_*,$ then the domains $D_{\alpha,r^-,r^+}$ are so-called \textit{elementary Reinhardt domains}.

Below we shall give a complete description of all proper holomorphic mappings between the domains $D_{\alpha,r^-_1,r^+_1}$ and
$D_{\beta,r^-_2,r^+_2}$ for arbitrary $\alpha,\beta\in(\R^2)_*$ and $0<r^+_i<\infty,\ -\infty<r^-_i<r^+_i,\ i=1,2.$ Similar problems were
studied in some papers. In \cite{shi1} and \cite{shi2} the problem of holomorphic equivalence of elementary Reinhardt domains was considered.
These results were partially extended by A.~Edigarian and W.~Zwonek. In the paper \cite{prop} the authors gave a characterization of proper
holomorphic mappings between elementary Reinhardt domains of the rational type.

\medskip Set
$\A(\rho^-,\rho^+):=\{z\in\C:\ \rho^-<|z|<\rho^+\}$ for $0<\rho^+,\ \rho^-<\rho^+$ and $\A_{\rho}:=\A(1/\rho,\rho),\ \rho>1.$ Moreover, put
\begin{align}\nonumber&D_{\gamma,r}:=\{z\in\C^2:\
1/r<|z_1||z_2|^{\gamma}<r\},\quad \gamma\in\R_*,\ r>1
\\\nonumber &D_{\gamma}:=\{z\in\C^2:\ |z_1||z_2|^{\gamma}<1\},\quad
\gamma\in\R_*
\\\nonumber&D_{\gamma}^*:=\{z\in\C^2:\
0<|z_1||z_2|^{\gamma}<1\},\quad\gamma\in\R_*
\end{align}

Note that if $\gamma$ is rational i.e. $\gamma=p/q$ for some relatively prime $p,q\in\Z,\ q>0,$ then $D_{\gamma,r}$ is biholomorphically
equivalent to $\A_{r^q}\times\C_*$ and $D_{\gamma}^*$ is  biholomorphically equivalent to $\D_*\times\C.$ Indeed, put
$$\psi(z_1,z_2):=(z_1^qz_2^p,z_1^mz_2^n)\quad \text{for}\quad (z_1,z_2)\in\C^2,$$
where $m,n\in\Z$ are such that $pm-qn=1.$ One can check that the mappings $\psi|_{D_{\gamma,r}}:D_{\gamma,r}\rightarrow \A_{r^q}\times\C_*$ and
$\psi|_{D_{\gamma}^*}:D_{\gamma}^*\rightarrow\D_*\times\C_*$ are biholomorphic.

Moreover, one may easily prove that $D_{\alpha,r^-,r^+}$ is algebraically equivalent to a domain of one of the following types:
\begin{enumerate} \item[(i)] If $r^->0$
\begin{enumerate}
\item[(a)]$\A_{\rho}\times\C,\quad\alpha_1\alpha_2=0,$ \item[(b)]$\A_{\rho}\times\C_*,\quad\alpha_1/\alpha_2\in\Q_*,$
\item[(c)]$D_{\gamma,\rho},\quad \gamma=\alpha_2/\alpha_1\in\NQ,$
\end{enumerate}
\item[(ii)]If $r^-=0$
\begin{enumerate}
\item[(a)]$\D_*\times\C,\quad\alpha_1\alpha_2=0,$ \item[(b)]$\D_*\times\C_*,\quad\alpha_1/\alpha_2\in\Q_*,$
\item[(c)]$D_{\gamma}^*,\quad\gamma=\alpha_2/\alpha_1\in\NQ,$
\end{enumerate}
\item[(iii)]If $r^-<0$
\begin{enumerate}
\item[(b)]$\D\times\C,\quad\alpha_1\alpha_2=0,$ \item[(a)]$D_{\gamma},\quad\gamma=\alpha_2/\alpha_1\neq0.$
\end{enumerate}
\end{enumerate}

Our main result is the following

\begin{theo}\label{nq}
(a) If $\alpha\in\NQ,$ then the set of proper holomorphic mappings between the domains $D_{\alpha,r}$ and $D_{\beta,R}$ is non-empty if and only
if
\begin{equation}\label{teza}\frac{\log R}{\log r}\in \Z + \beta \Z
\quad and \quad \alpha\frac{\log R}{\log r}\in \Z + \beta\Z.\end{equation}

(b) Let $\alpha,\beta\in\R\backslash\Q,$ $r,R>1$ be such that $\frac{\log R}{\log r}=k_1+l_1\beta$ and $\alpha\frac{\log R}{\log
r}=k_2+l_2\beta$ for some integers $k_i,l_i,\ i=1,2.$ Then any proper holomorphic mapping $f:D_{\alpha,r}\rightarrow D_{\beta,R}$ is given by
one of the following forms:
\begin{equation}\label{glowne}\begin{cases}f(z)=(az_1^{k_1}z_2^{k_2},bz_1^{l_1}z_2^{l_2})\quad or\\
f(z)=(az_1^{-k_1}z_2^{-k_2},bz_1^{-l_1}z_2^{-l_2})\end{cases}\ z=(z_1,z_2)\in D_{\alpha,r},\end{equation} where $a,b\in\C$ satisfy the relation
$|a||b|^{\beta}=1.$

Moreover, any of the mappings given by the formula (\ref{glowne}) is proper.
\end{theo}

Notice that in Theorem \ref{nq} (a) we do not demand $\beta$ to be irrational.

Using Theorem \ref{nq} we will easily obtain analogous results for the domains of the forms (ii) and (iii) of the irrational type.

\begin{theo}\label{el*}
Let $\alpha,\beta\in\NQ.$ The set of proper holomorphic mappings between the domains $D_{\alpha}^*$ and $D_{\beta}^*$ is non-empty if and only
if $\alpha=(k_2+\beta l_2)/(k_1+\beta l_1)$ for some $k_i,l_i\in\Z,\ i=1,2.$

Moreover, if $\alpha=(k_2+\beta l_2)/(k_1+\beta l_1),$ where $k_1+l_1\beta>0,$ then any proper holomorphic mapping $f:D_{\alpha}^*\rightarrow
D_{\beta}^*$ is of the form
\begin{equation}f(z_1,z_2)=(az_1^{k_1}z_2^{k_2},bz_1^{l_1}z_2^{l_2}),\quad (z_1,z_2)\in D_{
\alpha}^*, \end{equation} where $a,b\in\C$ satisfy the relation $|a||b|^{\beta}=1.$
\end{theo}

\begin{theo}\label{el}
Let $\alpha,\beta\in\NQ.$
\begin{enumerate}[(a)]
\item If $\alpha>0,\beta>0,$ then the set $\prop(D_{\alpha},D_{\beta})$ is non-empty if and only if $\alpha=p\beta$ for some $p\in\mathbb Q_{>0}.$ In this case all proper maps between $D_{\alpha}$ and $D_{\beta}$ are of the form \begin{equation}(z_1,z_2)\to(az_1^k,bz_2^l),\end{equation} where $a,b\in\mathbb C_*,\ |a||b|^{\alpha}=1,$ and $k,l$ are any positive integers satisfying the relation $p=\frac{l}{k}.$
\item If $\alpha<0,\beta<0,$ then the set $\prop(D_{\alpha},D_{\beta})$ is non-empty if and only if $\alpha=p_1+p_2\beta$ for some rational $p_1,p_2,\ p_2\neq0.$ In this case all proper maps between $D_{\alpha}$ and $D_{\beta}$ are of the form \begin{equation}(z_1,z_2)\to(az_1^{k_1}z_2^{k_2},bz_2^l),\end{equation} where $a,b\in\mathbb C_*,\ |a||b|^{\alpha}=1,$ and $k_1,k_2,l,\ k_1>0,$ are any integers satisfying the relations $p_1=-\frac{k_2}{k_1},\ p_2=\frac{l}{k_1}.$
\item If $\alpha\beta<0$ then there is no proper holomorphic mappings between $D_{\alpha}$ and $D_{\beta}.$
\end{enumerate}
\end{theo}

Next we prove the following

\begin{theo}\label{i}Let $\alpha,\beta\in(\R^2)_*,\ r^+_i>0,\ r^-_i<r^+_i,\ i=1,2.$
Assume that the sets $D_{\alpha,r^-_1,r^+_1},\ D_{\beta,r^-_2,r^+_2}$ are of the same type (either rational or irrational).

If there exists a proper holomorphic mapping $\psi:D_{\alpha,r^-_1,r^+_1}\rightarrow D_{\beta,r^-_2,r^+_2},$ then $r^-_1r^-_2>0$ or
$r^-_1=r^-_2=0.$
\end{theo}

In the case when the domains $D_{\alpha,r_1^-,r_1^+}$ and $D_{\alpha,r_2^-,r_2^+}$ are of different types we have the following result:

\begin{theo}\label{rozne}Let $\alpha,\beta\in(\R^2)_*,\ r^+_i>0,\ r_i^-<r^+_i,\
i=1,2.$ If the sets $D_{\alpha,r_1^-,r^+_1}$ and $D_{\beta,r^-_2,r^+_2}$ are of different types, then there is no proper holomorphic mapping
between $D_{\alpha,r_1^-,r^+_1}$ and $D_{\beta,r^-_2,r^+_2}.$
\end{theo}

Finally we discuss the rational case. As already mentioned the set of proper holomorphic mappings between elementary Reinhardt domains of the
rational type was described in \cite{prop}. Thus, in order to obtain the desired characterization, it suffices to prove the following three
theorems.
\begin{theo}\label{A} Let $r,R>1.$ If $R\neq r^m$ for any natural
$m,$ then the sets $\prop(\A_r\times\C,\A_R\times\C),$ $\prop(\A_r\times\C_*,\A_R\times\C)$ and $\prop(\A_r\times\C_*,\A_R\times\C_*)$ are
empty.

Moreover, for any $m\in\N:$

(a) $\prop(\A_r\times\C,\A_{r^m}\times\C)$ consists of the mappings of the form
$$\A_r\times\C\ni(z,w)\rightarrow(e^{i\theta}z^{\epsilon m},a_N(z)w^N+\ldots+a_0(z))\in\A_{r^m}\times\C,$$
where $\theta\in\R,\ N\in\N,\ \epsilon=\pm1$ and $a_0,\ldots, a_N\in\mathcal{O}(\A_r)$ are such that $|a_0(z)|+\ldots+|a_N(z)|>0,\ z\in \A_r.$

(b) $\prop(\A_r\times\C_*,\A_{r^m}\times\C)$ consists of the mappings of the form
$$\A_r\times\C_*\ni(z,w)\rightarrow(e^{i\theta}z^{\epsilon m},\frac{a_N(z)w^N+\ldots+a_0(z)}{w^k})\in\A_{r^m}\times\C,$$
where $\theta\in\R,\ k,N\in\N,\ 0<k<N,\ \epsilon=\pm1$ and $a_i\in\mathcal{O}(\A_r),\ i=1,\ldots,N$ satisfy the relations
$|a_0(z)|+\ldots+|a_{k-1}(z)|>0,\ |a_{k+1}(z)|+\ldots+|a_N(z)|>0,\ z\in \A_r.$

(c) $\prop(\A_r\times\C_*,\A_{r^m}\times\C_*)$ consists of the mappings of the form
$$\A_r\times\C_*\ni(z,w)\rightarrow(e^{i\theta}z^m,a(z)w^k)\in\A_{r^m}\times\C_*,$$
where $\epsilon=\pm1,\ \theta\in\R,\ k\in\N$ and $a\in\mathcal{O}(\A_r,\C_*).$
\end{theo}

\begin{theo}\label{A1}
There are no proper holomorphic mappings between the sets $\A_r\times\C$ and $\A_R\times\C_*$ for any $r,R>1.$
\end{theo}

\begin{theo}\label{A2}
(a) $\prop(\D_*\times\C,\D_*\times\C)$ consists of the mappings of the form
$$\D_*\times\C\ni(z,w)\rightarrow(e^{i\theta}z^m,a_N(z)w^N+\ldots+a_0(z))\in\D_*\times\C,$$
where $\theta\in\R,\ N\in\N,\ m\in\N$ and $a_0,\ldots, a_N\in\mathcal{O}(\D_*)$ are such that $|a_0(z)|+\ldots+|a_N(z)|>0,\ z\in\D_*.$

(b) $\prop(\D_*\times\C_*,\D_*\times\C)$ consists of the mappings of the form
$$\D_*\times\C_*\ni(z,w)\rightarrow(e^{i\theta}z^m,\frac{a_N(z)w^N+\ldots+a_0(z)}{w^k})\in\D_*\times\C,$$
where $\theta\in\R,\ m\in\N,\ k,N\in\N,\ 0<k<N$ and $a_i\in\mathcal{O}(\D_*),\ i=1,\ldots,N$ satisfy the relations
$|a_0(z)|+\ldots+|a_{k-1}(z)|>0,\ |a_{k+1}(z)|+\ldots+|a_N(z)|>0\ \text{for}\ z\in\D_*.$

(c) $\prop(\D_*\times\C_*,\D_*\times\C_*)$ consists of the mappings of the form
$$\D_*\times\C_*\ni(z,w)\rightarrow(e^{i\theta}z^m,a(z)w^k)\in\D_*\times\C_*,$$
where $\theta\in\R,\ k\in\N$ and $a\in\mathcal{O}(\D_*,\C_*).$

(d) The set $\prop(\D_*\times\C,\D_*\times\C_*)$ is empty.
\end{theo}

\section{Proofs}

The following result is probably known. However, we could not find it in the literature, so we present below our own proof.

\begin{lemma}\label{obserwacja}
Let $D\subset \C^n$ be a domain, $\alpha\in\NQ$ and let $f,g:D\rightarrow \C$ be holomorphic mappings satisfying the relation
$|f(z)|=|g(z)|^{\alpha},\ z\in D.$ Then, either $f=g=0$ on $D$ or there exists a holomorphic branch of logarithm $g$ i.e. a mapping
$\psi\in\mathcal{O}(D)$ such that $e^{\psi}=g$ on $D.$ In particular, there exists a $\theta\in\R$ such that $f=e^{i\theta+\alpha\psi}$ on $D.$
\end{lemma}

\begin{proof}
Comparing multiplicities of the roots of the functions $f$ and $g$ composed with affine mappings we may reduce our considerations to the case
when $f,g:D\rightarrow \C_*.$ Moreover, we may assume that $g(x')\in\R_{>0}$ for some $x'\in D.$

Obviously, there exists an $\eta\in\R$ such that the set $G_{\eta}:=\{z\in D:\ e^{i\eta}f(z)\in g(z)^{\alpha}\}$ is non-empty. Considering, if
necessary, a mapping $e^{i\eta}f$ instead of $f$ we may assume that $\eta=0.$

It is easy to see that $G_{0}$ in an open-closed subset of D, therefore, $G_0=D.$ Thus, there exists a holomorphic branch of $g^{\alpha}$ (also
denoted by $g^{\alpha}$) such that $g^{\alpha}(x')\in\R_{>0}$. It follows that there exist $g^t$ for any $t\in Q:=\{k+l\alpha:\ k,l\in\Z\}.$ Fix
a sequence $(t_m)_{m=1}^{\infty}\subset Q$ converging to $0$. In virtue of Montel's theorem, it is clear that $g^{t_m}\rightarrow 1$ locally
uniformly.

Put $\psi_m := \frac{g^{t_m}-1}{t_m}.$ Then $\lim_{m\rightarrow \infty}\psi_m(x')=\log g(x')$ and the sequence
$(\psi'_m)_{m=1}^{\infty}=(g^{t_m-1}g')_{m=1}^{\infty}$ is locally uniformly convergent with the limit $\frac{1}{g}g'.$ Thus the sequence
$(\psi_m)_{m=1}^{\infty}$ converges locally uniformly on D. Denote its limit by $\psi.$ By the Weierstrass theorem $\psi$ is holomorphic on D
and $\psi'=\lim_{m\rightarrow\infty}\psi'_m=\frac {1}{g}g'.$

Let $\widetilde{D}\subset D$ be any simply connected neighborhood of $x'.$ Let $\widetilde{\psi}$ be a holomorphic mapping on $\widetilde{D}$
such that $g|_{\widetilde{D}}=e^{\widetilde{\psi}}$ and $\widetilde{\psi}(x')=\log g(x').$ It is easy to see that $\widetilde{\psi}=\psi$ on
$\widetilde{D},$ therefore, by the identity principle we conclude that $g=e^{\psi}$ on $D.$
\end{proof}

\begin{lemma}\label{accumulation} Let $0<r^+_i,\ -\infty<r^-_i<r^+_i,\
i=1,2,\ \alpha,\beta\in\R.$ Let $(\lambda_n)_{n=1}^{\infty}\subset\A(r^-_1,r^+_1).$ Assume that the~mapping
$\phi:D_{(1,\alpha),r^-_1,r^+_1}\rightarrow D_{(1,\beta),r^-_2,r^+_2}$ is holomorphic and proper. Put
$$v(\lambda):=|\phi_1(\lambda,1)||\phi_2(\lambda,1)|^{\beta},\quad
\lambda\in \A(r^-_1,r^+_1).$$

If the sequence $(\lambda_n)_{n=1}^{\infty}$ has no accumulation points in $\A(r^-_1,r^+_1),$ then $(v(\lambda_n))_{n=1}^{\infty}$ has no
accumulation points in $\A(r^-_2,r^+_2).$
\end{lemma}

\begin{proof}
Assume that $v(\lambda_n)\rightarrow q.$ It suffices to show that $q\in\partial\A(r^-_2,r^+_2).$

Otherwise $q\in\A(r^-_2,r^+_2).$ Note that for any $\lambda\in\A(r^-_1,r^+_1)$ the function
\begin{equation}\label{foliacja} u_{\lambda}:\C\ni z\rightarrow|\phi_1(\lambda
e^{-\alpha z},e^{z})||\phi_2(\lambda e^{-\alpha z},e^{z})|^{\beta}\end{equation} is bounded and subharmonic, so $u_{\lambda}$ is constant.

Since $\phi$ is proper, the mapping $\C\ni z\rightarrow\phi_2(\lambda_n e^{-\alpha z},e^{z})\in\C$ is non-constant for any $n\in\N.$ Picard's
theorem implies that there is a sequence $(z_n)_{n=0}^{\infty}\subset\C$ such that $|\phi_2(\lambda_n e^{-\alpha z_n},e^{z_n})|^{\beta}=1.$
Obviously $u_{\lambda}(z)=u_{\lambda}(1)=v(\lambda),\ z\in\C$ and $v(\lambda_n)\rightarrow q,$ so $|\phi_1(\lambda_n e^{-\alpha
z_n},e^{z_n})|\rightarrow q.$ In particular, the set $\{\phi(\lambda_n e^{-\alpha z_n},e^{z_n}):\ n\in\N\}$ is relatively compact in
$D_{(1,\beta),r^-_2,r^+_2},$ however  the sequence $((\lambda_n e^{-\alpha z_n},e^{z_n}))_{n=1}^{\infty}$ does not have any accumulation points
in $D_{(1,\alpha),r_1^-,r_1^+};$ a contradiction.
\end{proof}

\begin{cor}\label{wlasciwosc} Let $\phi=(\phi_1,\phi_2):D_{\alpha,r}\rightarrow D_{\beta,R}$
be a proper holomorphic mapping and let $\alpha,\beta\in\R_{>0},\ r,R>1$. Put $v(\lambda):=|\phi_1(\lambda,1)||\phi_2(\lambda,1)|^{\beta},$
$\lambda\in\A_r.$ Then, either $$\lim_{|\lambda|\rightarrow 1/r}v(\lambda)=1/R,\ \lim_{|\lambda|\rightarrow r}v(\lambda)=R\ \quad{or}\
\lim_{|\lambda|\rightarrow 1/r}v(\lambda)=R,\ \lim_{|\lambda|\rightarrow r}v(\lambda)=1/R.$$
\end{cor}

\begin{lemma}\label{nql}
Let $\alpha\in\NQ,\ \beta\in\R,\ -\infty<r_i^-<r_i^+<\infty,\ 0<r_i^+,\ i=1,2,$ and let $\phi:D_{(1,\alpha),r^-_1,r^+_1}\rightarrow
D_{(1,\beta),r^-_2,r^+_2}$ be a~holomorphic mapping. Then for any $\lambda\in\A(r^-_1,r^+_1):$
\begin{align*}\phi(\{(z_1,z_2)\in&\C^2:\
|z_1||z_2|^{\alpha}=|\lambda|\})\subset\\\nonumber &\{(w_1,w_2)\in \C^2:\
|w_1||w_2|^{\beta}=|\phi_1(\lambda,1)||\phi_2(\lambda,1)|^{\beta}\}.\end{align*}
\end{lemma}

\begin{proof}
Note that for any $\lambda\in \A(r^-_1,r^+_1)$ the function
$$u:\C\ni z\rightarrow|\phi_1(\lambda e^{\alpha
z},e^{-z})||\phi_2(\lambda e^{\alpha z},e^{-z})|^{\beta}$$ is subharmonic and bounded. Hence $u$ is constant.

In virtue of Kronecker's theorem, the set $\{(|\lambda|e^{\alpha z},e^{-z}):z\in\C\}$ is dense in $\{(z_1,z_2)\in\C^2:\
|z_1||z_2|^{\alpha}=|\lambda|\}.$ Thus, there is $t\in\R$ such that
\begin{align*}\phi(\{(z_1,z_2)\in\C^2:\
|z_1||z_2|^{\alpha}=|\lambda|\})\subset\nonumber \{(w_1,w_2)\in \C^2:\ |w_1||w_2|^{\beta}=t\}.\end{align*}

It is easy to see that $t=|\phi_1(\lambda,1)||\phi_2(\lambda,1)|^{\beta}.$
\end{proof}

\begin{proof}[Proof of Theorem \ref{nq} (a)] Let $\phi:D_{\alpha,r}\rightarrow
D_{\beta,R}$ be a proper holomorphic mapping. Put $v(\lambda):=|\phi_1(\lambda,1)||\phi_2(\lambda,1)|^{\beta},\ \lambda\in\A_r.$ Obviously,
$\log v$ is a harmonic function. Applying Corollary \ref{wlasciwosc} and Hadamard's theorem we infer that $v$ is one of the following forms
$$v(\lambda)=|\lambda|^{\frac{\log R}{\log r}},\quad \lambda\in\A_r\quad \text{or}\quad
v(\lambda)=|\lambda|^{-\frac{\log R}{\log r}},\quad \lambda\in\A_r.$$ From this and Lemma \ref{nql} we easily conclude that there is
$\epsilon=\pm 1$ such that
\begin{equation}\label{had}|\phi_1(z)||\phi_2(z)|^{\beta}=|z_1|^{\epsilon\frac{\log R}{\log
r}}|z_2|^{\epsilon\alpha\frac{\log R}{\log r}},\quad z\in D_{\alpha,r}.\end{equation}

Let $z_2=1,\ z_1=z\in\A_r,\ \psi_i(z):=\phi_i(z,1),\ i=1,2.$ Then
$$\log(\psi_1(z)\overline{\psi}_1(z))+\beta\log(\psi_2(z)\overline{\psi}_2(z))=\epsilon\frac{\log
R}{\log r}\log(z\overline{z}).$$  Differentiating with respect to $z$ we get
\begin{equation}\label{der}\frac{\psi'_1(z)}{\psi_1(z)}+\beta\frac{\psi'_2(z)}{\psi_2(z)}=\epsilon\frac{\log R}{\log
r}\cdot\frac{1}{z},\quad z\in\A_r. \end{equation}

It follows that
$$\ind(\psi_1\circ\gamma;0)+\beta\ind(\psi_2\circ\gamma;0)=\epsilon\frac{\log R}{\log
r},$$ where $\gamma$ is the unit circle. Hence $\frac{\log R}{\log r}\in \Z + \beta \Z.$ The same argument with respect to the second variable
shows that $\alpha\frac{\log R}{\log r}\in \Z + \beta\Z.$

To prove the converse, assume that the conditions (\ref{teza}) are fulfilled, i.e.
$$\frac{\log R}{\log r}=k_1+l_1\beta,\quad \alpha\frac{\log R}{\log r}=k_2+l_2\beta,$$
where $k_i,l_i\in\Z,\ i=1,2.$ Define $\phi_1(z):=z_1^{k_1} z_2^{k_2},\ \phi_2(z):=z_1^{l_1} z_2^{l_2}\quad \text{for}\quad z=(z_1,z_2)\in\C^2
,\quad \phi:=(\phi_1,\phi_2).$ Observe that $\phi|_{D_{\alpha,r}}\in\prop(D_{\alpha,r};D_{\beta,R}).$ Indeed, it is easy to check that
\begin{equation}\label{wlas}|\phi_1(z)||\phi_2(z)|^{\beta}=|z_1|^{\log
R/\log r}|z_2|^{\alpha\log R/\log r},\quad (z_1,z_2)\in D_{\alpha,r},\end{equation} so
$\phi|_{D_{\alpha,r}}\in\mathcal{O}(D_{\alpha,r},D_{\beta,R}).$ Since $k_1l_2\neq k_2l_1,$ $\phi$ is a proper holomorphic mapping from
$(\C_{\ast})^2$ into itself (see \cite{wlasciwe}, Theorem 2.1). Now we immediately conclude from (\ref{wlas}) that $\phi|_{D_{\alpha,r}}$ is a
proper holomorphic mapping between $D_{\alpha,r}$ and $D_{\beta,R}$.
\end{proof}

Lemma \ref{obserwacja} and (\ref{had}) lead to the following

\begin{cor}\label{warunek konieczny} Let $\alpha,\beta\in\NQ,$ and
$\phi\in\prop(D_{\alpha,r},D_{\beta,R}).$ Assume that $\frac{\log R}{\log r}=k_1+l_1\beta$ and $\alpha\frac{\log R}{\log r}=k_2+l_2\beta$ for
some $k_i,l_i\in\Z,\ i=1,2.$  Then there are $\theta\in\R,\ \psi\in\mathcal{O}(D_{\alpha,r})$ and $\epsilon\in\{1,-1\}$ such that

$$\phi(z)=(z_1^{\epsilon k_1}z_2^{\epsilon
k_2}e^{i\theta}e^{-\beta\psi(z)},z_1^{\epsilon l_1}z_2^{\epsilon l_2}e^{\psi(z)}),\quad z\in D_{\alpha,r}.$$
\end{cor}

\begin{remark}\label{pm}  We
may always assume that $\epsilon$ in Corollary \ref{warunek konieczny} is equal to $1$ (if necessary instead of $\phi$ we may consider the
mapping $\phi\circ h,$ where $h\in\aut(D_{\alpha,r})$ is given by the formula $h(z_1,z_2):=(z_1^{-1},z_2^{-1})$).
\end{remark}

To prove Theorem \ref{nq} we need the following notation. Put $X_{\alpha, r}:=\{z\in\C^2:\ -\log r< Re z_1+\alpha Re z_2<\log r\}$ and
$\Pi(z_1,z_2):=(e^{z_1},e^{z_2})$ for $ (z_1,z_2) \in \C^2.$ Note that $(X_{\alpha,r},\Pi)$ is the universal covering of $D_{\alpha,r}.$
Moreover, it is clear that $X_{\alpha,r}$ is simply connected.

We get the following lemma

\begin{lemma}\label{podniesienie}
Let $\alpha,\beta\in\NQ,\ r,R>1$ and assume that $\frac{\log R}{\log r}=k_1+l_1\beta,\ \alpha\frac{\log R}{\log r}=k_2+l_2\beta,$ where
$k_i,l_i\in\Z,\ i=1,2.$ Let $f:D_{\alpha,r}\rightarrow D_{\beta,R}$ be a proper holomorphic mapping. Then every continuous lifting of the
mapping $f\circ\Pi:X_{\alpha,r}\rightarrow D_{\beta,R}$ is proper and holomorphic.
\end{lemma}

\begin{proof}
In virtue of Corollary \ref{warunek konieczny} and Remark \ref{pm} we may assume that the mapping $f$ is given by the formula
$f(z)=(z_1^{k_1}z_2^{k_2}e^{-\beta\psi(z)+i\theta},z_1^{l_1}z_2^{l_2}e^{\psi(z)}),\ (z\in D_{\alpha,r}),$ where $\theta\in \R$ and
$\psi\in\mathcal{O}(D_{\alpha,r}).$

Let $\widetilde{f}$ be any continuous lifting of $f\circ\Pi:X_{\alpha,r}\rightarrow D_{\beta,R},$ that is $\widetilde{f}:X_{\alpha,r}\rightarrow
X_{\beta,R}$ and $f\circ\Pi=\Pi\circ\widetilde{f}.$ It is obvious that $\widetilde{f}$ is holomorphic. Then by the identity principle
\begin{equation}\label{podniesienie e}\begin{cases}\widetilde{f}_1(z)=k_1z_1+k_2z_2-\beta\psi(e^{z_1},e^{z_2})+i\theta+2\mu_1\pi
i,\\
\widetilde{f}_2(z)=l_1z_1+l_2z_2+\psi(e^{z_1},e^{z_2})+2\mu_2\pi i\end{cases} z\in X_{\alpha,r}\end{equation} for some $\mu_i\in\Z,\ i=1,2.$

Suppose that $\widetilde{f}$ is not proper, i.e. there is a~sequence $(z^m)_{m=1}^{\infty}\subset X_{\alpha,r},\ z^m=(z_1^m,z_2^m),\ m\in\N$
without any accumulation points in $X_{\alpha,r}$ such that $(\widetilde{f}(z^m))_{m=1}^{\infty}$ is convergent in $X_{\beta,R}.$ Put
$y_0:=\lim_{m\rightarrow\infty}\widetilde{f}(z^m)\in X_{\beta,R}.$

Obviously, $f(\Pi(z_1^m,z_2^m))=\Pi(\widetilde{f}(z_1^m,z_2^m))\rightarrow \Pi(y_0).$ Since $f$ is proper, the set $\{\Pi(z^m):\ m\geq 1\}$ is
relatively compact in $D_{\alpha,r}.$ Thus we may assume that the sequence $(\Pi(z^m))_{m=1}^{\infty}$ is convergent in $D_{\alpha,r}.$ Denote
its limit by $w_0:=\lim_{m\rightarrow\infty}\Pi(z^m)\in D_{\alpha,r}.$ From (\ref{podniesienie e}) we deduce that the sequences
$(k_1z_1^m+k_2z_2^m)_{m=1}^{\infty}$ and $(l_1z_1^m+l_2z_2^m)_{m=1}^{\infty}$ are convergent in $\C^2.$ Thus $(z^m)_{m=1}^{\infty}$ is also
convergent.

Put $z_0:=\lim_{m\rightarrow\infty}z^m.$ Now it suffices to observe that $\Pi(z_0)=w_0\in D_{\alpha,r},$ so $z_0\in X_{\alpha,r};$ a
contradiction
\end{proof}

Now we are able to give a description of the set of proper holomorphic mappings between the domains $D_{\alpha, r}$ and $D_{\beta,R}$ of the
irrational type.

\begin{proof}[Proof of Theorem \ref{nq} (b)]

Let $f\in\prop(D_{\alpha,r},D_{\beta,R})$. In virtue of Corollary \ref{warunek konieczny} and Remark \ref{pm} we may assume that
$$f(z)=(z_1^{k_1}z_2^{k_2}e^{-\beta\psi(z)+i\theta},z_1^{l_1}z_2^{l_2}e^{\psi(z)}),\quad z=(z_1,z_2)\in D_{\alpha,r}$$
for some $\theta\in\R$ and $\psi\in\mathcal{O}(D_{\alpha,r}).$ Our aim is to show that the mapping $\psi$ is constant.

To simplify notation for $\gamma\in\R$ put
$$\Lambda_{\gamma}:\C^2\ni(z_1,z_2)\rightarrow(z_1+\gamma
z_2,z_2)\in\C^2.$$ It is clear that $\Lambda_{\gamma}(X_{\gamma,\rho})=S_{\rho}\times\C,\ \rho>1,$ where $S_{\rho}:=\{z\in\C:\ -\log r<Re z<\log
r\}.$ Moreover, the mapping $\Lambda_{\gamma}$ is biholomorphic and the inverse is given by $\Lambda_{\gamma}^{-1}=\Lambda_{-\gamma}.$

Note that the mapping $\widetilde{f}:X_{\alpha,r}\rightarrow X_{\beta, R}$ given by
$$\widetilde{f}(z)=(k_1z_1+k_2z_2-\beta\psi(e^{z_1},e^{z_2})+i\theta,l_1z_1+l_2z_2+\psi(e^{z_1},e^{z_2}))$$
is a lifting of $f\circ\Pi.$ Thus Lemma \ref{podniesienie} implies that $\widetilde{f}$ is proper and holomorphic.

Put $H:=(H_1,H_2):=\Lambda_{\beta}\circ\widetilde{f}\circ\Lambda_{\alpha}^{-1}:S_r\times\C\rightarrow S_R\times\C.$ Obviously, the mapping $H$
is proper and holomorphic.

Applying the relations $\frac{\log R}{\log r}=k_1+l_1\beta,\quad \alpha\frac{\log R}{\log r}=k_2+l_2\beta$ we see that
\begin{equation}\label{h}H(z)=(z_1(k_1+\beta
l_1)+i\theta,l_1z_1+z_2(l_2-l_1\alpha)+\psi(e^{z_1-\alpha z_2},e^{z_2})),\quad z\in S_r\times \C.\end{equation}

From this we conclude that for any  $z_1\in S_r$ the mapping $\C\ni z\rightarrow H_2(z_1,z)\in \C$ is proper and holomorphic. Consequently, due
to the form of proper holomorphic self-mappings of $\C,$ we deduce that there is a polynomial $p=p_{z_1}\in \mathcal{P}(\C)$ such that
$H_2(z_1,z)=p(z).$ Therefore, the polynomial $q(z):=q_{z_1}(z):=p(z)-l_1z_1-z(l_2-l_1\alpha)$ satisfies the equation
\begin{equation}\psi(e^{z_1}e^{-\alpha z},e^z)=q(z),\quad z\in\C.\end{equation}

Notice that $\{( e^{z_1} e^{- \alpha 2\pi i m}, e^{2 \pi i m}):\ m\in\N\}$ is a relatively compact subset of $D_{\alpha,r}$ and the sequence
$\{q(2\pi i m)\}_{m=1}^{\infty}$ is bounded. Thus the polynomial $q$ is constant.

Put $c(z_1):=\psi(e^{z_1-\alpha z_2},e^{z_2}),\ z_1\in S_r.$ Let us fix any $1<\rho<R$ and take a~constant $M=M(\rho)>0$ such that $|c(x)|<M$
for every $x\in[-\log\rho,\log\rho].$

Let $\lambda\in\rho \mathbb{D}\setminus\frac{1}{\rho}\mathbb{D}$ be arbitrary. Note that for any $z_2\in\C$ we have $|\psi(|\lambda|e^{-\alpha
z_2},e^{z_2})|=|c(\log|\lambda|)|<M.$ Applying Kronecker's theorem we infer that the set $\{(|\lambda|e^{-\alpha z},e^z):\ z\in \C\}$ is dense
in $\{(z_1,z_2)\in\C^2:\ |z_1||z_2|^{\alpha}=|\lambda|\}.$ Consequently $\psi|_{D_{\alpha,\rho}}$ is bounded.

Now it suffices to repeat the proof of Lemma 2.7.1 of \cite{metryki} in order to show that every bounded holomorphic mapping on
$D_{\alpha,\rho}$ (in particular $\psi$) is constant.

On the other hand, we have already mentioned in the proof of Theorem \ref{nq} (a), that any mapping given by the formula (\ref{glowne}) is
proper.
\end{proof}

\begin{proof}[Proof of Theorems \ref{el*} and \ref{el}]
We prove simultaneously both cases. Let $f:D_{\alpha}\rightarrow D_{\beta}$ (respectively, $f:D_{\alpha}^*\rightarrow D_{\beta}^*$) be a proper
holomorphic function. We aim at reducing the situation to that of Theorem \ref{nq}. Take any $r>1.$

From Lemma \ref{nql} we see that for any $t\in[0,1)$ ($t\in(0,1)$) there is an $s(t)\in[0,1)$ ($s(t)\in(0,1)$) such that
$$f(\{(z_1,z_2)\in\C^2:\ |z_1||z_2|^{\alpha}=t\})\subset\{(w_1,w_2)\in\C^2:\
|w_1||w_2|^{\beta}=s(t)\}.$$ Note that $s(|\lambda|)=|f_1(\lambda,1)||f_2(\lambda,1)|^{\beta}$ and the function $v$ given by
$v:\D\ni\lambda\rightarrow s(|\lambda|)\in\ [0,1]$ (respectively $v:D_*\ni\lambda\rightarrow s(|\lambda|)\in\ [0,1]$) is radial and subharmonic
on $\D$ (in the second case we may remove singularity at $0$). The maximum principle applied to the function $v$ implies that $s$ is increasing.

In particular, there is an $R>1$ such that the restriction $f|_{D_{(1,\alpha),1/r^2,1}}:D_{(1,\alpha),1/r^2,1}\rightarrow D_{(1,\beta),1/R^2,1}$
is proper. For $\rho>1$ put $\widetilde{\Lambda}_{\rho}:\C^2\ni(z_1,z_2)\rightarrow(\rho z_1,z_2)\in\C^2$ and define
$\psi:=\widetilde{\Lambda}_R\circ f\circ\widetilde{\Lambda}_r^{-1}|_{D_{\alpha,r}}.$ Note that $\psi:D_{\alpha,r}\rightarrow D_{\beta,R}$ is a
proper holomorphic mapping. Applying Theorem \ref{nq} we find that that $\frac{\log R}{\log r}=k_1+l_1\beta,\ \alpha\frac{\log R}{\log
r}=k_2+l_2\beta$ and $\psi(z_1,z_2)=(a z_1^{\epsilon k_1}z_2^{\epsilon k_2},b z_1^{\epsilon l_1}z_2^{\epsilon l_2})$ for some $k_i,l_i\in\Z,\
i=1,2,\ \epsilon=\pm 1$ and $a,b\in\C$ satisfying the equation $|a||b|^{\beta}=1.$ Obviously \begin{equation}\label{er}\alpha=\frac{k_2+l_2\beta}{k_1+l_1\beta}\end{equation} and by
the identity principle we obtain that
\begin{equation}f(z_1,z_2)=(ar^{\epsilon l_1\beta}z_1^{\epsilon k_1}z_2^{\epsilon k_2},
br^{-\epsilon l_1}z_1^{\epsilon l_1}z_2^{\epsilon l_2}),\quad (z_1,z_2)\in D_{\alpha}\ ((z_1,z_2)\in D_{\alpha}^*).\end{equation}

If we consider the case $f:D_{\alpha}^*\rightarrow D_{\beta}^*,$ then it suffices to notice that
$|f_1(z)||f_2(z)|^{\beta}=(|z_1||z_2|^{\alpha})^{\epsilon k_1+\epsilon l_1\beta},\ z=(z_1,z_2)\in D_{\alpha}^*,$ hence $\epsilon (k_1+
l_1\beta)>0.$

Now let us focus our attention on the remaining case, i.e. the situation when $f:D_{\alpha}\rightarrow D_{\beta}.$ 

It is clear that if $\alpha\beta<0$ then the set $\prop(D_{\alpha},D_{\beta})$ is empty.

Assume that $\alpha,\beta>0.$ Considerations done at the beginning of the proof show that $f$ preserves the axes $(\mathbb C\times\{0\})\cup(\{0\}\times\mathbb C).$ Therefore $f(z_1,z_2)=(az_1^{k_1},bz_2^{l_2}),\ (z_1,z_2)\in D_{\alpha},$ or $f(z_1,z_2)=(az_2^{l_1}b,z_1^{k_2}),\ (z_1,z_2)\in D_{\alpha},$ where $k_i,l_i\geq0,\ i=1,2.$ A direct computation shows that $f$ cannot be of the second form (otherwise by (\ref{er}) we would find that $\alpha\in\mathbb Q$). From this piece of information one can easily get (a).

Similarly, if $\alpha,\beta<0$ we state that the mapping $f$ is of the form $f(z_1,z_2)=(az_1^{k_1}z_2^{k_2},bz_2^{l_2}),\ (z_1,z_2)\in D_{\alpha},\ k_1\geq0.$ As before, using this piece of information one can easily finish the proof.

From this piece of information we easily get the required formulas.

On the other hand, one can check that any of the mappings given in Theorem \ref{el} is proper (since $\alpha$ is irrational,  
$k_1l_2-k_2l_1\neq0$).
\end{proof}

\begin{lemma}\label{li}
Let $r^+>0,\ r^-<r^+,\ t\in\R.$ Suppose that the function $v:\A(r^-,r^+)\rightarrow[-\infty,t)$ is subharmonic, radial (i.e.
$v(|\lambda|)=v(\lambda),\ \lambda\in \A(r^-,r^+)$) and harmonic on the set $\{z\in \A(r^-,r^+):\ v(z)\neq-\infty\}.$ Then there exist
$a,b\in\R$ such that
$$v(\lambda)=a\log|\lambda|+b,\quad \lambda\in \A(r^-,r^+).$$\end{lemma}

\begin{proof}It suffices to observe that since $v$ is radial,
$\A(r^-,r^+)\setminus\{0\}\subset\{z\in \A(r^-,r^+):\ v(z)\neq-\infty\}$ (and next one may proceed standardly, i.e. solve an easy differential
equation).
\end{proof}

\begin{proof}[Proof of Theorem \ref{i}]
First, let us consider the case when $D_{\alpha,r^-_1,r^+_1}$ and $D_{\beta,r^-_2,r^+_2}$ are of the irrational type. Then we may assume that
$\alpha=(1,\alpha_1)$ for some $\alpha_1\in\NQ.$ Let
$$v:\A(r^-_1,r^+_1)\ni\lambda\rightarrow\log|\psi_1(\lambda,1)|^{\beta_1}|\psi_2(\lambda,1)|^{\beta_2}\in\R.$$
By Lemma \ref{nql} we get that $\psi(\{(z_1,z_2)\in\C^2:\ |z_1||z_2|^{\alpha}=|\lambda|\}\subset\{(w_1,w_2)\in\C^2:\
|w_1|^{\beta_1}|w_2|^{\beta_2}=e^{v(\lambda)}\}.$ Therefore, the function $v$ is radial. Observe moreover that $v$ is subharmonic on
$\A(r^-_1,r^+_1)$ and harmonic on the set $\{\lambda\in\A(r^-_1,r^+_1):\ v(\lambda)>-\infty\}.$ Since $\psi$ is surjective, we conclude that
\begin{equation}\label{obraz}v(\A(r^-_1,r^+_1))=\begin{cases}(\log r^-_2,\log r^+_2),\quad \text{if}\ r^-_2\geq0,\\
[-\infty,\log r^+_2),\quad \text{if}\ r^-_2<0,\end{cases}\end{equation} (we put $\log0:=-\infty$). However, by Lemma \ref{li} the function $v$
must be of the form $v(\lambda)=a\log|\lambda| + b,\ \lambda\in\A(r^-_1,r^+_1)$ for some $a,b\in\R,$ which easily finishes the proof in this
case.

Now suppose that $D_{\alpha,r^-_1,r^+_1}$ and $D_{\beta,r^-_2,r^+_2}$ are of the rational type; without loss of generality we may assume that
$\beta=(p,q)\in\Z^2$ and $\alpha=(1,\alpha_1)$ for some $\alpha_1\in\Q.$ Applying Lemma \ref{accumulation} one can see that the mapping
$$\A(r^-_1,r^+_1)\ni\lambda\rightarrow\psi_1(\lambda,1)^p\psi_2(\lambda,1)^q\in\A(r^-_2,r^+_2)$$
is proper. Hence this case follows directly from the form of the set of proper holomorphic mappings between $\A(r^-_1,r^+_1)$ and
$\A(r^-_2,r^+_2).$
\end{proof}

\begin{proof}[Proof of Theorem \ref{rozne}]
Assume that $D_{\alpha,r^-_1,r^+_1}$ is of the rational type and $D_{\beta,r^-_2,r^+_2}$ is of the irrational type; without loss of generality
$\alpha=(1,p/q)$ for some $p,q\in\Z$ and $\beta=(1,\beta_2)$ for some $\beta_2\in\NQ$.

Suppose that $\psi:D_{\alpha,r^-_1,r^+_1}\rightarrow D_{\beta,r^-_2,r^+_2}$ is a proper holomorphic mapping. Note that for any
$\lambda\in\A(r^-_1,r^+_1)$ the mapping
\begin{equation}\label{lambda} u_{\lambda}:\C_*\ni z\rightarrow|\psi_1(\lambda z^p,z^{-q})||\psi_2(\lambda z^p,z^{-q})|^{\beta_2} \end{equation} is constant.
Fix $\lambda_0$ and $c\neq0$ such that $u_{\lambda_0}\equiv c.$ One can see that $\C_*\ni z\rightarrow\psi_i(\lambda_0z^p,z^{-q})\in\C_*$ is a
proper holomorphic self-mapping of $\C_*,\ i=1,2.$ Therefore, there are $a_i\in\C_*$ and $\mu_i\in\Z_*,\ i=1,2$ such that
$\psi_i(\lambda_0z^p,z^{-q})=a_iz^{\mu_i},\ z\in\C_*,\ i=1,2.$ Applying (\ref{lambda}) it is clear that
$|a_1||a_2|^{\beta_2}|z|^{\mu_1+\mu_2\beta_2}=c,\ z\in\C_*.$ In particular, $\beta_2\in\Q;$ a contradiction.

Now, suppose that there exists a proper holomorphic mapping $\psi:D_{\beta,r^-_2,r^+_2}\\ \rightarrow D_{\alpha,r^-_1,r^+_1}.$ Put
$u(\lambda):=|\psi_1(\lambda,1)||\psi_2(\lambda,1)|^{\beta_2}$ for $\lambda\in\A(r^-_2,r^+_2).$

Applying Lemmas \ref{accumulation} and \ref{nql} we obtain that the function $u$ satisfies assumptions of Lemma \ref{li}. Thus, there are
$a,b\in\R$ such that $\log u(\lambda)=a\log|\lambda|+b,\ \lambda\in\A(r^-_2,r^+_2).$ In particular, the function $u$ is either strictly
increasing or strictly decreasing. Take any $\rho^-_2$, $\rho^+_2$ such that $\rho^-_2>\max\{0,r^-_2\},$ $\rho_2^+<r^+_2,$
 $\rho_2^-<\rho_2^+.$ Put
$\rho^-_1:=\min\{u(\rho^-_2),u(\rho^+_2)\},$ $\rho^+_1:=\max\{u(\rho^-_2),u(\rho^+_2)\}.$ Then
$$\psi|_{D_{\beta,\rho^-_2,\rho^+_2}}:D_{\beta,\rho^-_2,\rho^+_2}\rightarrow
D_{(1,\alpha),\rho^-_1,\rho^+_1}$$ is obviously a proper holomorphic mapping. In virtue of Theorem \ref{nq} (a) we see that there are
$k_i,\l_i\in\Z,\ i=1,2$ such that $\beta=(k_1+l_1\alpha)/(k_2+l_2\alpha).$ In particular, $\beta\in\Q;$ a contradiction.
\end{proof}

\begin{lemma}\label{A0}Let $A,B\subset\C^n$ be domains and assume that $B$ is bounded.

(a) The mapping $f:A\times\C_*\rightarrow B\times\C$ is proper and holomorphic if and only if there are $m\in\prop(A,B),\ k\in\N,\ 0<k<N,\
N\in\N,\ a_i\in\mathcal{O}(A),\ i=1,\ldots,N,\ |a_0(z)|+\ldots+|a_{k-1}(z)|>0,\ |a_{k+1}(z)|+\ldots+|a_N(z)|>0,\ z\in A$ satisfying the relation
$$f(z,w)=\left(m(z),\frac{a_N(z)w^N+\ldots+a_0(z)}{w^k}\right),\quad (z,w)\in A\times\C_*.$$

(b) The mapping $f:A\times\C\rightarrow B\times\C$ is proper and holomorphic if and only if there are $a_0,\ldots, a_N\in\mathcal{O}(A),\
N\in\N,\ |a_0(z)|+\ldots+|a_N(z)|>0,\ z\in A$ and there is a proper holomorphic mapping $m:A\rightarrow B$ such that
$$f(z,w)=(m(z),a_N(z)w^N+\ldots+a_0(z)),\quad (z,w)\in A\times\C.$$

(c) The mapping $f:A\times\C_*\rightarrow B\times\C$ is proper and holomorphic if and only if there are $m\in\prop(A,B),\
a\in\mathcal{O}(A,\C_*)$ and $k\in\N$ such that
$$f(z,w)=(m(z),a(z)w^k),\quad (z,w)\in A\times\C_*.$$

(d) There is no proper holomorphic mappings between $A\times\C$ and $B\times\C_*$.
\end{lemma}

\begin{proof}
First of all, notice that for any $z\in A$ the mapping $w\rightarrow f_1(z,w)\in\C^n$ is bounded on $\C$ (or $\C_*$), so it is constant.

\textit{(a)} Observe that $\C_*\ni w\rightarrow f_2(z,w)\in\C$ is a proper mapping for any $z\in A.$ Thus, for any $z\in A$ there is a
polynomial $p(z,\cdot),\ p(z,0)\neq 0,$ and a natural $k(z)$ such that
\begin{equation}\label{r}\phi_2(z,w)=\frac{p(z,w)}{w^{k(z)}},\quad
(z,w)\in A\times\C_*.\end{equation} One can see that there is a $k$ such that $k=k(z),\ z\in A$ (use Rouche's theorem). Consequently
$p\in\mathcal{O}(A\times\C_*)$.

Fix any domain $A'\subset\subset A$ and put
$$A_{\mu}:=\{z\in\overline{A'} :\ \frac{\partial^{\mu} p}{\partial
w^{\mu}}(z,w)=0 \ \text{for any}\ w\in\C\}.$$ The above considerations imply that $\bigcup_{\mu=1}^{\infty}A_{\mu}=\overline{A'}.$ Applying
Baire's theorem we find that there exits $N\in\N$ such that $A_{N}$ does not have empty interior. By the identity principle $A_{N}=A.$

Thus, there are holomorphic mappings $a_0,\ldots,a_N:A\rightarrow\C$ such that $p(z,w)=a_N(z)w^N+\ldots+a_1(z)w+a_0(z)$ for $(z,w)\in
A\times\C,$ i.e.
\begin{equation}f_2(z,w)=\frac{a_N(z)w^N+\ldots+a_1(z)w+a_0(z)}{w^k},\quad
(z,w)\in A\times\C.\end{equation} By properness of $f_2(z,\cdot)$ we conclude that $0<k<N,\ |a_N(z)|+\ldots +|a_{k+1}(z)|>0$ and
$|a_{k-1}(z)|+\ldots+|a_0(z)|>0$ for any $z\in A.$

Put $m(z):=f_1(z,1),\ z\in A.$ We claim that m is proper.

Indeed, take any sequence $(z_n)_{n=1}^{\infty}$ and assume that it does not have any accumulation points in $A.$ Without loss of generality we
may assume that $a_0(z_n)\neq0$ for any $n\in\N$ (if required we may replace $a_0$ with $a_1$ etc.). Then there exists a sequence
$(w_n)_{n=1}^{\infty}\subset\C_*$ such that $a_N(z_n)w_n^N+\ldots+a_1(z_n)w_n+a_0(z_n)=0$ for any $n\in\N.$ Since $f(z_n,w_n)=(m(z_n),0),$ it is
obvious that $(m(z_n))_{n=1}^{\infty}$ has no accumulation points in $B.$

Conversely one can check that every mapping $f$ defined in this way is proper.

\textit{(b)} It is easy to see that $\C\ni w\rightarrow f_2(z,w)\in\C$ is the proper holomorphic mapping for any $z\in A.$ From the form of
proper holomorphic self-mappings we conclude that for every $z\in A$ the mapping $f_2(z,\cdot)$ is a complex polynomial. Now we proceed exactly
as in the proof of (a).

\textit{(c)} We proceed similarly as in the proof of (a) and (b).

\textit{(d)} Suppose that $f:A\times\C\rightarrow B\times\C_*$ is a proper holomorphic function. Fix $z\in A.$ Then the mapping $\C\ni
w\rightarrow f_2(z,w)\in\C_*$ is proper.

Take $\psi\in\mathcal{O}(\C)$ such that $f_2(1,\cdot)=\exp\circ\psi.$ Observe that $\psi$ is a proper holomorphic self-mapping of the complex
plane, hence $\psi$ is a polynomial. From these we easily get a contradiction.
\end{proof}

\begin{proof}[Proof of Theorems \ref{A}, \ref{A1} and \ref{A2}]
It is a direct consequence of Lemma \ref{A0}.
\end{proof}
\medskip
Finally, I take this opportunity to express deep gratitude to professor W\l odzimierz Zwonek for introducing me to the subject and numerous
remarks.

\end{document}